\documentclass{elsart}
\usepackage{latexsym,amssymb}
\usepackage{flafter}
\usepackage{graphicx}
\usepackage{algorithm}
\usepackage{algorithmic}

\newtheorem{theorem}{Theorem}[section]

\begin{document}
\runauthor{Jing Meng, Pei-yong Zhu, Hou-Biao Li}
\begin{frontmatter}
\title{A SQMRCGstab Algorithm for Families of Shifted Linear Systems\thanksref{b}}
\thanks[b]{\small Supported by the National Natural Science Foundation of China (11026085, 11101071, 10671134, 11271001, 51175443)
and the Fundamental Research Funds for China Scholarship Council.}

\author{Jing Meng\corauthref{a}},
\corauth[a]{Corresponding author.}
\ead{lihoubiao0189@163.com (H-B Li), zpy6940@sina.com (P-Y Zhu)}
\ead{mengmeng-erni@163.com}
\author{Pei-Yong Zhu},
\author{Hou-Biao Li}
\address{School of Mathematical Sciences, University of Electronic Science and Technology of
China, Chengdu, 611731, P. R. China}

\begin{abstract}
This study is mainly focused on iterative solutions to shifted linear systems arising from a Quantum Chromodynamics (QCD) problem.
To solve such system efficiently, we explore a kind of shifted QMRCGstab (SQMRCGstab) methods,
which is derived by extending the quasi-minimum residual to the shifted BiCGstab. The shifted QMRCGstab method takes advantage of the shifted structure,
so that the number of matrix-vector products and the number of inner products are the same as a single linear system.
Moreover, the SQMRCGstab achieves a smoothing of the residual compared to the shifted BiCGstab, and is more competitive than the MS-QMRIDR(s) and the shifted BiCGstab on the QCD problem.
Numerical examples show also the efficiency of the method when one applies it to the real problems.
\end{abstract}


\end{frontmatter}

\section{Introduction}

Quantum Chromodynamics (QCD) is generally accepted to be the fundamental physical theory of strong interactions among the quarks as constituents of matter. To explore some physical observables in QCD,
it is important to discretize the Dirac operator. It could respect the chiral symmetry which such observables depend on.
To study QCD at nonzero baryon density, the more powerful overlap Dirac operator was extended to include a quark chemical potential recently in \cite{Bloch1}.
The most important and challenging part in the overlap operator is to compute the sign function of a complex matrix,
which is Hermitian at zero baryon density, but is non-Hermitian at non-zero chemical potential.

To compute the sign function of the large non-Hermitian spare matrices efficiently, Jacques C.R. Blocha, etc. \cite{Bloch2,Bloch3,Bloch4} used Multi-shift methods to
approximate the sign function $f$ by a rational function $g$
\[f(t)\approx g(t)=\sum_{i=1}^{s}\frac{\omega_{i}}{t-\sigma_{i}}.\]
Then it implies that
\[f(A)b=\sum_{i=1}^{s}\omega_{i}x^{i},\]
in which the $x^{i}, i=1,2,\ldots s,$ are solutions of the $s$ linear systems
\begin{equation}\label{eq:1}
 (A-\sigma_{i} I)x^{i}=b, ~~(i=1,2,\ldots s)
\end{equation}
where $I$ is the identity matrix, the large nonsingular and non-Hermitian matrix $A\in \mathbb{C}^{n\times n}$, the parameter $\sigma_{i}\in \mathbb{C}$ and $\sigma_{i} \notin \lambda(A)$ ($\lambda(A)$ is the set of all eigenvalues of $A$), $1\leq i\leq s$, $b\in \mathbb{C}^{n}$.

In the present paper, we pay close attention to the simultaneous solutions of
the $s$ linear system for several tabulated values of $\sigma_{i},1\leq i\leq s$. For simply, we denote $\sigma =\{-\sigma_{i}|1\leq i\leq s\}$ and $x^{\sigma}=\{x^{i}|1\leq i\leq s\}$, then Eqs. \ref{eq:1} can be written as the following format
 \begin{equation}\label{eq:2}
 (A+\sigma I)x^{\sigma}=b,
\end{equation}
which is called a shifted linear system. The linear system $Ax=b$ will be termed as seed linear system. Sequences of such shifted linear systems arise not only in QCD \cite{Bloch2,Bloch3,Bloch4,Bloch5,Sakurai}, but also in various fields, e.g., in trust-region and regularization techniques for nonlinear least squares and other optimization
problems in control theory \cite{Datta}, as well as in the application of implicit methods for the numerical
solution of partial differential equations (PDEs) \cite{Gallopoulos}. However, the QCD application is the main motivation for the present study.

To solve (\ref{eq:1}), Krylov subspace techniques are most faithfully. Since they rely on a shift-invariance property, which allows to obtain approximation iterates for all parameter values by only constructing one approximation subspace.

Recently, many of Krylov subspace methods for shifted systems were proposed. For example, the shifted CG \cite{Van}, COCG \cite{Takayama} and CGCR \cite{Sogabe} were considered for the case of systems with Hermitian coefficient matrices. The shifted restarted FOM and restarted GMRES were powerful solver for non-Hermitian systems. In addition, the restarted GMRES, which forced the shifted system residual to be colinear to the seed system residual, modified the GMRES iteration for the shift system \cite{Frommer1}. Unfortunately, these methods may result in a numerically unstable process, so that after a few restarts numerical results become useless, for detail, see \cite{Frommer1,Simoncini}.

Since the quasi-minimum residual is not shift invariant, the QMR, TFQMR and MINRE methods do not define their iterates by a Petrov-Galerkin condition. However, these methods allow to save the matrix-vector multiplication for the shifted system. It is due to the fact that these methods construct a basis for $K_{m}(A,b)$ via the Lanczos process, and this basis is invariant under shifts, see \cite{Freund1} and also \cite{Jegerlehner}.
In 2003, A. Frommer \cite{Frommer} proposed the shift BiCGStab ($\ell$) and showed that for a positive real matrix $A$ and a positive shift $\sigma$, the proposed method was a well-smoothed variant of BiCG. Taken $\ell=1$, the shifted BiCGstab method is obtained.

The shifted BiCGstab  \cite{Jegerlehner,Frommer}, among the shifted algorithms, is a particularly efficient method for quark propagator calculation. However, its convergence curve is not smoothed (see Section 5). In order to eliminate that erratic convergence, we derive an alternative approach (SQMRCGstab), which is applied the quasi-minimun residual to the shifted BiCGstab and illustrate its smoothed convergence by means of numerical experiments in this paper.
The proposed algorithm also makes use of the special shifted structure, and for any family of shifted systems, the number of matrix-vector products and the number of inner products are the same as a single linear system.
Moreover, we also compare the SQMRCGstab method with the Multi-shift QMRIDR(s)\cite{MartinB} method. Multi-shift QMRIDR(s) (MS-QMRIDR(s)) is the Quasi-Minimal Residual variant of the IDR(s)\cite{Martin} for solving shifted systems.
The numerical experiments show that the SQMRCGstab is more competitive than the MS-QMRIDR(s) and the shifted BiCGstab on the QCD problem.

The rest of this paper is organized as follows. In Section 2, we briefly review collinear residuals principle. The shifted BiCGstab algorithm is shortly summarized in Section 3 . Our algorithm details are described in Section 4. Section 5 presents numerical examples, which stem form real problems such as QCD, to illustrate our results.

\section{Collinear Residuals}

Given a seed system  $Ax=b$ and an initial vectors $x_{0}$, we consider a biorthogonal Krylov subspace method for the iterative solutions.
Let $V_{m}$ and $W_{m}$, which are built by the Lanczos biorthonalization algorithm, be a pair of biorthognal bases for the two subspaces
\[
K_{m}(A,r_{0})=span\{r_{0},Ar_{0},\ldots,A^{m-1}r_{0}\}
\]
and
\[
K_{m}(A^{H},r_{0})=span\{r_{0},A^{H}r_{0},\ldots,(A^{H})^{m-1}r_{0}\},
\]
where $r_{0}=b-Ax_{0}$. The first basis vector $v_{1}$ of $V_{m}$ is $r_{0}/\|r_{0}\|$, and then the following relation holds,
 \begin{equation}\label{eq:3.1}
 AV_{m}=V_{m}T_{m}+t_{m+1,m}v_{m+1}e_{m}^{H},
 \end{equation}
in which the matrix $T_{m}$ is the projection of $A$ obtained from an oblique projection process onto $K_{m}(A,r_{0})$ and orthogonally to $K_{m}(A^{H},r_{0})$.
Next, let us consider the shift system (\ref{eq:2}). The (\ref{eq:3.1}) will be transformed into the following form by shifting
\[
(A+\sigma I)V_{m}=V_{m}(T_{m}+\sigma I)+t_{m+1,m}v_{m+1}e_{m}^{H}.
\]
The above character is the well-know shift-invariant property of Krylov subspace. That is,
\[K_{m}(A,r_{0})=K_{m}(A+\sigma I,r_{0})=K_{m}(A-\sigma_{i} I,r_{0}),~~1\leq i\leq s.\]
Therefore, it allows to obtain approximation iterates for all parameter values by only constructing one approximation subspace,
which allows to save the matrix-vector multiplication for solving the $s$ linear systems.

A Krylov subspace method produces iterative solutions $x_{m}$ for which the residuals $r_{m}=b-Ax_{m}$ are in the Krylov space $K_{m}(A,r_{0})$. As a consequence, the residual $r_{m}$ can be written as $p_{m}(A)r_{0}$, where $p_{m}$ is a polynomial of degree $\leq m-1$ with $p_{m}(0)=1$.

Similarly, for the shifted system, any vector $x_{m}^{\sigma}$ forming the affine Krylov subspace $K_{m}(A+\sigma I,r_{0}^{\sigma})$ can be represented as
$x_{m}^{\sigma}=x_{0}^{\sigma}+q_{m-1}^{\sigma}(A+\sigma I)r_{0}^{\sigma}\in x_{0}^{\sigma}+K_{m}(A+\sigma I,r_{0}^{\sigma})$,
where $q_{m-1}^{\sigma}$ is a polynomial of degree $\leq m-1$. The corresponding residual $r_{m}^{\sigma}=b-(A+\sigma I)x_{m}^{\sigma}$ satisfies
\[
r_{m}^{\sigma}=r_{0}^{\sigma}-(A+\sigma I)q_{m-1}^{\sigma}(A+\sigma I)r_{0}^{\sigma}=p_{m}^{\sigma}(A+\sigma I)r_{0}^{\sigma},
\]
where $p_{m}^{\sigma}(t+\sigma )=1-(t+\sigma )q_{m-1}^{\sigma}(t+\sigma )$ is a polynomial of degree $\leq m$ with $p_{m}^{\sigma}(0)=1$.

For convenience, we assume that any Krylov subspace method is started with a starting guess $x_{0}=0$ or $x_{0}^{\sigma}=0$, thus $r_{0}^{\sigma}=r_{0}=b$.
Krylov subspace methods for the shifted system exploit the following result, from which we learn what is a collinear residual idea. For the following proof we refer to \cite{Frommer}.
\begin{theorem}(\cite{Frommer}).
Let $\Gamma_{1}\subseteq\Gamma_{2}\subseteq\ldots\subseteq\Gamma_{m}$ be a sequence of nested subspaces of $\mathbb{C}^{n}$ (i.e., test spaces) such that $\Gamma_{i}$ has dimension $i$ and $\Gamma_{i}\cap(K_{i+1}(A,b))^{\bot}={0}, \ i=1,\ldots,m.$ Let $x_{i}\in K_{i}(A,b)$ be an approximation to the solution of $Ax=b$ defined via the following Petrov-Galerkin condition for the residual $r_{i}=b-Ax_{i}=p_{i}(A)b$:
$
r_{i}\bot\Gamma_{i},  i=1,\ldots,m.
$
Similarly, let $x_{i}^{\sigma}\in K_{i}(A+\sigma I,b)=K_{i}(A,b)$ be the approximation to the solution of $(A+\sigma I)x^{\sigma}=b$ with residual $r_{i}^{\sigma}=b-(A+\sigma I)x_{i}^{\sigma}=p_{i}^{\sigma}(A+\sigma I)b$, again satisfying
$
r_{i}^{\sigma}\bot\Gamma_{i},  i=1,\ldots,m.
$
Then $r_{i}$ and $r_{i}^{\sigma}$ are collinear, i.e.
\[
 r_{i}^{\sigma}=(1/c_{i}^{\sigma})r_{i}, \   c_{i}^{\sigma}\in \mathbb{C}.
\]
\end{theorem}

To solve the shifted system (\ref{eq:2}), we take the above collinear residual approach as follows,
\begin{equation}\label{eq:2.3}
 r_{m}^{\sigma}=(1/c_{m}^{\sigma})r_{m}, \   c_{m}^{\sigma}\in \mathbb{C}.
\end{equation}
The (\ref{eq:2.3}) is equivalent to
$p_{m}^{\sigma}(A+\sigma I)b=(1/c_{m}^{\sigma})p_{m}(A)b.$
By Comparing coefficients, the follow identity is obtained
\[
p_{m}^{\sigma}(t+\sigma )=(1/c_{m}^{\sigma})p_{m}(t).
\]

Since $p_{m}^{\sigma}(0)=1$, the equation $c_{m}^{\sigma}=p_{m}(-\sigma)$ is obtained. If the parameter $\sigma>0$, then $c_{m}^{\sigma}=p_{m}(-\sigma)>1$, see \cite{Frommer,Vann}.

\section{The Shifted BiCGstab Algorithm}

In this section, we briefly recall the derivation of the shifted BiCGstab method, which is similar to the one described in \cite{Jegerlehner} or \cite{Frommer}.
At first, we shortly describe the algorithm of BiCG when it is applied to a seed system $Ax=b$ (see Algorithm 1).

{\bf Algorithm 1: The BiCG Algorithm (\cite{Saad}).}
\begin{enumerate}
  \item Computer $r_{0}=b-Ax_{0}$. Choose $r_{0}^{*}$ such that $(r_{0},r_{0}^{*})\neq 0$.
  \item Set, $u_{0}=r_{0},  u_{0}^{*}=r_{0}^{*}$.
  \item {\bf for} $m=0,1,\ldots$, until convergence {\bf do}
  \begin{itemize}
    \item $a_{m}=(r_{m},r_{m}^{*})/(Au_{m}^{*},u_{m}^{*})$
    \item $x_{m+1}=x_{m}+\alpha_{m}u_{m}$
    \item $r_{m+1}=r_{m}-\alpha_{m}Au_{m}$
    \item $r_{m+1}^{*}=r_{m}^{*}-\alpha_{m}A^{\top}u_{m}^{*}$
    \item $\beta_{m}=(r_{m+1},r_{m+1}^{*})/(r_{m},r_{m}^{*})$
    \item $u_{m+1}=r_{m+1}+\beta_{m}u_{m}$
    \item $u_{m+1}^{*}=r_{m+1}^{*}+\beta_{m}u_{m}^{*}$
  \end{itemize}
  \item {\bf end for}
\end{enumerate}

From the 6th step and the 9th one in Algorithm 1, one knows that the following recurrence relation holds,
\begin{equation}\label{eq:4}
\begin{array}{lll}
r_{m+1}&=& (1+\frac{\beta_{m-1}}{\alpha_{m-1}}\alpha_{m}-\alpha_{m}A)r_{m}-\frac{\beta_{m-1}}{\alpha_{m-1}}\alpha_{m}r_{m-1}\\
&=&-\alpha_{m}Ar_{m}+ (1+\frac{\beta_{m-1}}{\alpha_{m-1}}\alpha_{m})r_{m}-\frac{\beta_{m-1}}{\alpha_{m-1}}\alpha_{m}r_{m-1}.
\end{array}
\end{equation}
Similarly, applying the BiCG method to the shifted system $(A+\sigma I)x^{\sigma}=b$, we may consider updating the residual $r_{m+1}^{\sigma}$ of the shifted system with the following three-term recurrence relation
\begin{equation}\label{eq:5}
r_{m+1}^{\sigma}= (1+\frac{\beta_{m-1}^{\sigma}}{\alpha_{m-1}^{\sigma}}\alpha_{m}^{\sigma}-\alpha_{m}^{\sigma}(A+\sigma I))r_{m}^{\sigma}-\frac{\beta_{m-1}^{\sigma}}{\alpha_{m-1}^{\sigma}}\alpha_{m}^{\sigma}r_{m-1}^{\sigma}.
\end{equation}
where $\alpha_{m}^{\sigma}, \beta_{m}^{\sigma}$ are the analogous to the coefficient $\alpha_{m}, \beta_{m}$ in BiCG algorithm.

To obtain the computational formula for $r_{m+1}^{\sigma}$, three parameters $\alpha_{m}^{\sigma},\beta_{m}^{\sigma}$ and $c_{m+1}^{\sigma}$ must be completely fixed. Hence, we next give computational formulas for the three values. Substituting (\ref{eq:2.3}) into (\ref{eq:5}), we obtain
\begin{equation}\label{eq:5.1}
\begin{array}{lll}
r_{m+1}&=&-\alpha_{m}^{\sigma}(\frac{c_{m+1}^{\sigma}}{c_{m}^{\sigma}})Ar_{m}+ (1+\frac{\beta_{m-1}^{\sigma}}{\alpha_{m-1}^{\sigma}}\alpha_{m}^{\sigma}-\alpha_{m}^{\sigma}\sigma I)\frac{c_{m+1}^{\sigma}}{c_{m}^{\sigma}}r_{m}\\
&&-\frac{\beta_{m-1}^{\sigma}\alpha_{m}^{\sigma}c_{m+1}^{\sigma}}{\alpha_{m-1}^{\sigma}c_{m-1}^{\sigma}}r_{m-1}.
\end{array}
\end{equation}
Then by comparing the coefficients with (\ref{eq:4}), the computational formulas for the three parameters may be obtained
\[
c_{m+1}^{\sigma}=(1+\alpha_{m}\sigma)c_{m}^{\sigma}+\frac{\alpha_{m}\beta_{m}}{\alpha_{m-1}}(c_{m-1}^{\sigma}-c_{m}^{\sigma}),
\]
\[
\alpha_{m}^{\sigma}=\alpha_{m}(\frac{c_{m}^{\sigma}}{c_{m+1}^{\sigma}}),\]
\[
\beta_{m}^{\sigma}=(\frac{c_{m-1}^{\sigma}}{c_{m}^{\sigma}})^2\beta_{m}.
\]

It is well known that, by combining the BiCG method and the minimal residual (MR) idea together, the BiCGstab algorithm \cite{Vorst} was derived. During each step of the BiCGstab method for solving the shifted linear systems, the following relations hold,
\begin{equation}\label{eq:6}
s_{m}^{\sigma}=r_{m}^{\sigma}-\alpha_{m}^{\sigma}(A+\sigma I)u_{m}^{\sigma},   \             \
r_{m+1}^{\sigma}=s_{m}^{\sigma}-\chi_{m}^{\sigma}(A+\sigma I)s_{m}^{\sigma},
\end{equation}

where $s_{m}^{\sigma}$, $\chi_{n}^{\sigma}$ are the analogous to $s_{m}$ and $\omega_{m}$ respectively in BiCGstab \cite{Saad}, $u_{m}^{\sigma}$ is the search direction, and $\chi_{n}^{\sigma}$ is chosen by a local steepest descent principle. Therefore, $r_{m}^{\sigma}$ and  $s_{m}^{\sigma}$ satisfy the following equations,
\[
r_{m}^{\sigma}=\psi_{m}^{\sigma}(A+\sigma I)p_{m}^{\sigma}(A+\sigma I)r_{0}^{\sigma},
\]
\[
s_{m}^{\sigma}=\psi_{m-1}^{\sigma}(A+\sigma I)p_{m}^{\sigma}(A+\sigma I)r_{0}^{\sigma},
\]
in which $p_{m}^{\sigma}$ is the degree m residual polynomial of the m-th step of the BiCG method for the shifted linear system and $\psi_{m}^{\sigma}$ is the MR polynomial which is defined recursively at each step with the goal of stabilizing or smoothing convergence behavior.
Since $p_{m}^{\sigma}(A+\sigma I)b=(1/c_{m}^{\sigma})p_{m}(A)b$, we obtain the following form,
\[
r_{m}^{\sigma}=\psi_{m}^{\sigma}(A+\sigma I)(1/c_{m}^{\sigma})p_{m}(A)b,
\]
\[
s_{m}^{\sigma}=\psi_{m-1}^{\sigma}(A+\sigma I)(1/c_{m}^{\sigma})p_{m}(A)b.
\]

To obtain the computational formulas for $r_{m}^{\sigma}$ and $s_{m}^{\sigma}$, the parameters $\chi^{\sigma}_{m}$ and the shifted MR polynomial $\psi_{m}^{\sigma}$ are required. In the following, we will show how to calculate the parameters $\chi^{\sigma}_{m}$ and the shifted MR polynomial $\psi_{m}^{\sigma}$, for more details, see \cite{Jegerlehner}.

By ignoring the scalar coefficients at first, the MR polynomial $\psi_{m}(t)$ is defined by the simple recurrence
$\psi_{m+1}(t)=(1-\chi_{m}t)\psi_{m}(t)$ in the BiCGstab method. Then the polynomial is given directly as a product of its linear factor $\psi_{m+1}(t)=\prod_{i=1}^{m}(1-\chi_{i}t)$.
To calculate the shifted MR polynomial, we assume a linear factor
 $(1-(t+\sigma)\chi^{\sigma})=c(1-t\chi)$,
resulting in
$\chi^{\sigma}=\frac{\chi}{1+\sigma\chi},\ c=\frac{1}{1+\sigma\chi}$.
The shifted polynomial is therefore given by
\begin{equation}\label{eq:7}
\begin{array}{lll}
\psi_{m}^{\sigma}(t+\sigma)&=&\prod^{m}_{i=1}(1-(t+\sigma)\chi_{i}^{\sigma})\\
&=&\prod^{m}_{i=1}(1-(t+\sigma)\frac{\chi_{i}}{1+\sigma\chi_{i}})\\
&=&\prod^{m}_{i=1}\frac{1}{1+\sigma\chi_{i}}(1-\chi_{i}t)\\
&=&\rho_{m}^{\sigma}\psi_{m}(t),\\
\end{array}
\end{equation}
\begin{equation}\label{eq:8}
\rho_{m}^{\sigma}=\prod_{i=1}^{m}\frac{1}{1+\sigma\chi_{i}}.
\end{equation}
Thus using the polynomials (\ref{eq:7}) and (\ref{eq:8}), we can generate the parameters as follows
\[
\chi_{m}^{\sigma}=\frac{\chi_{m}}{1+\sigma\chi_{m}},
\]
\[
\rho_{m+1}^{\sigma}=\frac{\rho_{m}^{\sigma}}{1+\sigma\chi_{m}}.
\]

In addition, $x_{m}^{\sigma}$ is computed from $x_{m-1}^{\sigma}$ by  \[x_{m}^{\sigma}={x}_{m-1}^{\sigma}+\alpha_{m}^{\sigma}u_{m}^{\sigma}+\chi_{m}^{\sigma} s_{m}^{\sigma}.\]
Therefore, the $u_{m}^{\sigma}$ is essentially required. The problem is that the update of $u_{m}^{\sigma}$ required the calculation of $(A+\sigma I)u_{m}^{\sigma}$, which will increase $s$ additional matrix-vectors multiplication for Eq.(\ref{eq:1}). In order to eliminate these additional matrix-vector products, the following identify is used in the shifted BiCGstab (see Algorithm 2)
\[(A+\sigma I)u_{m}^{\sigma}=\frac{1}{\alpha_{m-1}^{\sigma}}(r_{m-1}^{\sigma}-s_{m}^{\sigma}).\]

The procedure described in this Section 3 leads to Algorithm 2.

{\bf Algorithm 2: A Shifted BiCGstab Algorithm (\cite{Jegerlehner,Frommer}).}
\begin{enumerate}
  \item Choose $\sigma\in S= \{-\sigma_{i},1\leq i\leq s\} $,
  \item Computer $r_{0}=b-Ax_{0}$. Choose $r_{0}^{*}$ such that $(r_{0},r_{0}^{*})\neq 0$.
  \item set $x_{0}=x_{0}^{\sigma}=\textbf{0}, r_{0}=\tilde{r}_{0}=r_{0}^{\sigma}=b$, $c_{0}^{\sigma}=c_{-1}^{\sigma}=\rho_{0}=\alpha_{0}=\chi_{0}=q_{0}^{\sigma}=1$,
$u_{0}=v_{0}=d_{0}=\textbf{0}$. 
  \item {\bf for} $m=0,1,\ldots$, until convergence {\bf do}
  \begin{itemize}
    \item  $\rho_{m}=(\tilde{r}_{0},r_{m-1})$,
    \item $\beta_{m}=\frac{\rho_{m}\alpha_{m-1}}{\rho_{m-1}\chi_{m-1}}$
    \item $u_{m}=r_{m-1}+\beta_{m}(u_{m-1}-\chi_{m-1}v_{m-1})$
    \item $v_{m}=Au_{m}$
    \item $\alpha_{m}=\rho_{m}/(\tilde{r}_{0},v_{m})$
    \item $s_{m}=r_{m-1}-\alpha_{m}v_{m}$
    \item $c_{m}^{\sigma}=(1+\alpha_{m}\sigma)c_{m-1}^{\sigma}+\frac{\alpha_{m}\beta_{m}}{\alpha_{m-1}}(c_{m-2}^{\sigma}-c_{m-1}^{\sigma})$
    \item $\alpha_{m}^{\sigma}=\alpha_{m}(\frac{c_{m-1}^{\sigma}}{c_{m}^{\sigma}}),$
    \item $\beta_{m}^{\sigma}=(\frac{c_{m-2}^{\sigma}}{c_{m-1}^{\sigma}})^2\beta_{m}. $
    \item $u_{m}^{\sigma}=r_{m-1}^{\sigma}+\beta_{m}^{\sigma}(u_{m-1}^{\sigma}-\chi_{m-1}^{\sigma}v_{m-1}^{\sigma})$
    \item $s_{m}^{\sigma}=c_{m}^{\sigma}q_{m-1}^{\sigma}s_{m}$
    \item $v_{m}^{\sigma}=\frac{1}{\alpha_{m-1}^{\sigma}}(r_{m-1}^{\sigma}-s_{m}^{\sigma})$
    \item compute $\chi_{m}^{\sigma},q_{m}^{\sigma}$ and update $r_{m},r_{m}^{\sigma}$
    \item $t_{m}=As_{m}, \ \chi_{m}=(s_{m},t_{m})/(t_{m},t_{m})$
    \item $\chi_{m}^{\sigma}=\frac{\chi_{m}}{1+\chi_{m}}$
    \item $q_{m}^{\sigma}=q_{m-1}^{\sigma}/(1+\chi_{m})$
    \item $r_{m}= s_{m}-\chi_{m}t_{m}$
    \item $r_{m}^{\sigma}=c_{m}^{\sigma}q_{m}^{\sigma}r_{m}$
    \item $x_{m}^{\sigma}={x}_{m-1}^{\sigma}+\alpha_{m}^{\sigma}u_{m}^{\sigma}+\chi_{m}^{\sigma} s_{m}^{\sigma}$
   \item If $x_{m}^{\sigma}$ is accurate enough, then quit
\end{itemize}
  \item {\bf end for}
\end{enumerate}

The shifted BiCGstab is a particularly efficient method for quark propagator calculation. However, its convergence curve is not smoothed (see Section 5). In older to eliminate that erratic convergence, we derive a method which is applied the quasi-minimun residual to the shifted BiCGstab in next section.

\section{ The SQMRCGstab Algorithm}

The algorithm proposed in this section is inspired by the QMRCGstab method \cite{Chan}, which was combined the features of BiCGstab and quasi-minimization principle. Note that if applying the quasi-minimization principle to the shifted BiCGSTAB method, one may obtain SQMRCGstab method. Equally, the collinear residual approach is extended to the QMRCGstab method for solving shifted linear systems.

Assuming that the vectors $r_{m}^{\sigma}$, $u_{m}^{\sigma}$ and $s_{m}^{\sigma}$ are generated by the shifted BiCGstab method, we choose $x_{m}^{\sigma}$ by quasi-minimizing the residual over their span.

Let
$Y_{m}^{\sigma}=\{y_{0}^{\sigma},y_{1}^{\sigma},\ldots,y_{m}^{\sigma}\}$, where $y_{2i-1}^{\sigma}=u_{i}^{\sigma}$, $i=1,2,\ldots,[(m+1)/2]$\footnote{[~~] is the integer function} and $y_{2i}^{\sigma}=s_{i}^{\sigma}$, $i=1,2,\ldots,[m/2]$. Similarly,
$W_{m}^{\sigma}=\{w_{0}^{\sigma},w_{1}^{\sigma},\ldots,w_{m}^{\sigma}\}$ with $w_{2i}^{\sigma}=r_{i}^{\sigma}$, $i=0,1,\ldots,[(m+1)/2]$ and $w_{2i-1}^{\sigma}=s_{i}^{\sigma}$, $i=1,2,\ldots,[m/2]$.
We also define $\{\delta_{1},\delta_{2},\ldots,\delta_{m}\}$, as $\delta_{2i}^{\sigma}=\chi_{i}^{\sigma}$ for $i=1,2,\ldots,[(m+1)/2]$ and $\delta_{2i-1}^{\sigma}=\alpha_{i}^{\sigma}$, $i=1,2,\ldots,[(m+1)/2]$. In this case, Eq.(\ref{eq:6}) may be written as
\[
(A+\sigma I)y_{i}^{\sigma}=(w_{i-1}^{\sigma}-w_{i}^{\sigma})\delta_{i}^{-1}, \\i=1,\ldots,m.
\]
By the definitions of $W_{m+1}^{\sigma}$ and $Y_{m}^{\sigma}$, the following relation is obvious that
\[
(A+\sigma I)Y_{m}^{\sigma}= W_{m+1}^{\sigma} \triangle_{m+1}^{\sigma}
\]
in which $\triangle_{m+1}^{\sigma}$is a $(m+1)\times m$ bidiagonal matrix, i.e.,
\[
\left(\begin{array}{cccccc}
    \delta_{1}^{-1} & 0 & \ldots & \ldots & 0 \nonumber  \\
    -\delta_{1}^{-1} & \delta_{2}^{-1} &  &   & \vdots \nonumber   \\
    0 & -\delta_{2}^{-1} & \delta_{3}^{-1} & \ldots&   \nonumber  \\
    \vdots &   &\ddots & \ddots &\vdots \nonumber  \\
    \vdots &   &   & -\delta_{m-1}^{-1} & \delta_{m}^{-1}  \nonumber \\
    0 & \ldots &   &   & -\delta_{m}^{-1}\nonumber  \\
\end{array}\right).
\]

It can be easily checked that span\{$Y_{m}^{\sigma}$\}=span$\{W_{m}^{\sigma}\}$=$K_{m}(A+\sigma I,r_{0}^{\sigma})$, where span$\{Y_{m}^{\sigma}\}$ is generated by the shifted BiCGSTAB method.

Next we use the quasi-minimization principle method to find an approximation to the solution (\ref{eq:2}) over span$\{Y_{m}^{\sigma}\}$. Apparently, the approximate solution $x_{m}^{\sigma}$ can be given by
\[
x_{m}^{\sigma}=x_{0}^{\sigma}+Y_{m}^{\sigma} z, \  \ z\in \mathbb{C}^{k}.
\]
Hence, the residual can be written as
\[
r_{m}^{\sigma}=r_{0}^{\sigma}+(A+\sigma I) Y_{m}^{\sigma}z=r_{0}^{\sigma}- W_{m+1}^{\sigma} \triangle_{m+1}^{\sigma}z.
\]
Since the first vector of $ W_{m+1}^{\sigma}$ is $r_{0}^{\sigma}=b$, it follows that
\[
r_{m}^{\sigma}=W_{m+1}^{\sigma}(e_{1}-  \triangle_{m+1}^{\sigma}z),
\]
where $e_{1}$ is the first vector of the canonical basis. In order to make the columns of $W_{m+1}^{\sigma}$ to be unit norm, we use a $(m+1)\times(m+1)$ scaling matrix $\Sigma_{m+1}=diag(\theta_{1},\theta_{2},\ldots,\theta_{m+1})$ with $\theta_{i}=\|w_{i}^{\sigma}\|$. Then
\begin{equation}\label{eq:9}
r_{m}^{\sigma}=W_{m+1}^{\sigma}\Sigma_{m+1}^{-1}(\theta_{1}e_{1}- H_{m+1}^{\sigma}z)
\end{equation}
with $H_{m+1}^{\sigma}=\Sigma_{m+1}\triangle_{m+1}^{\sigma}$.

In order to minimize the residual norm over the Krylov subspace, the 2-norm of the right-hand side of Eq.(\ref{eq:9}) would have to be minimized, but this is not practical since the columns of $W_{m+1}^{\sigma}\sum_{m+1}^{-1}$ are not orthonormal as in Arnoldi. However, $\|\theta_{1}e_{1}- H_{m+1}^{\sigma}z\|$ can be minimized over $z$, as was done for the QMR algorithm.

In this paper, the least squares minimization of $\|\theta_{1}e_{1}- H_{m+1}^{\sigma}z\|$ is solved using QR decomposition of $H_{m+1}^{\sigma}$. Since $H_{m+1}^{\sigma}$ is lower bidiagonal, this is done by means of Givens rotations, and only the rotation of the previous step is needed, for detail, see \cite{Saad}.

Finally, the SQMRCGstab algorithm is summarized as follows (see Algorithm 3), in which the Givens rotations used in the QR decomposition are given explicitly.

{\bf Algorithm 3: A SQMRCGstab Algorithm.}
\begin{enumerate}
  \item set $x_{0}=x_{0}^{\sigma}=\textbf{0}, r_{0}=\tilde{r}_{0}=r_{0}^{\sigma}=b$, $c_{0}^{\sigma}=c_{-1}^{\sigma}=\rho_{0}=\alpha_{0}=\chi_{0}=q_{0}^{\sigma}=1$,
        $\tau=\|r\|,\theta_{0}=\eta\mathcal{0}=0$,$u_{0}=v_{0}=d_{0}=\textbf{0}$.
  \item {\bf for} $m=0,1,\ldots$, until convergence {\bf do}
  \begin{itemize}
  \item $\rho_{m}=(\tilde{r}_{0},r_{m-1})$,
  \item $\beta_{m}=\frac{\rho_{m}\alpha_{m-1}}{\rho_{m-1}\chi_{m-1}}$ 
  \item $u_{m}=r_{m-1}+\beta_{m}(u_{m-1}-\chi_{m-1}v_{m-1})$
  \item $v_{m}=Au_{m}$
  \item $\alpha_{m}=\rho_{m}/(\tilde{r}_{0},v_{m})$
  \item $s_{m}=r_{m-1}-\alpha_{m}v_{m}$
  \item $c_{m}^{\sigma}=(1+\alpha_{m}\sigma)c_{m-1}^{\sigma}+\frac{\alpha_{m}\beta_{m}}{\alpha_{m-1}}(c_{m-2}^{\sigma}-c_{m-1}^{\sigma})$
  \item $\alpha_{m}^{\sigma}=\alpha_{m}(\frac{c_{m-1}^{\sigma}}{c_{m}^{\sigma}}),$
  \item $\beta_{m}^{\sigma}=(\frac{c_{m-2}^{\sigma}}{c_{m-1}^{\sigma}})^2\beta_{m}. $
\item  $u_{m}^{\sigma}=r_{m-1}^{\sigma}+\beta_{m}^{\sigma}(u_{m-1}^{\sigma}-\chi_{m-1}^{\sigma}v_{m-1}^{\sigma})$
\item $s_{m}^{\sigma}=c_{m}^{\sigma}q_{m-1}^{\sigma}s_{m}$
\item $v_{m}^{\sigma}=\frac{1}{\alpha_{m-1}^{\sigma}}(r_{m-1}^{\sigma}-s_{m}^{\sigma})$
\item first quasi-minimization and update iterate
\item $\tilde{\theta_{m}}=\|s_{m}^{\sigma}\|/\tau; \ \zeta=1/\sqrt{1+\tilde{\theta_{m}}^{2}}$
\item $\tilde{\tau}=\tau\tilde{\theta_{m}}\zeta; \   \tilde{\eta_{m}}=\zeta^{2}\alpha_{m}^{\sigma}$
\item $\tilde{d}_{m}^{\sigma}=u_{m}^{\sigma}+\frac{\theta_{m-1}^{2}\eta_{m-1}}{\alpha_{m}^{\sigma}}d_{m-1}$
\item $\tilde{x}_{m}^{\sigma}= x_{m-1}^{\sigma}+\tilde{\eta_{m}}\tilde{d}_{m}$
\item compute $\chi_{m}^{\sigma},q_{m}^{\sigma}$ and update $r_{m}$
\item $t_{m}=As_{m}, \ \chi_{m}=(s_{m},t_{m})/(t_{m},t_{m})$
\item $\chi_{m}^{\sigma}=\frac{\chi_{m}}{1+\chi_{m}}$
\item $q_{m}^{\sigma}=q_{m-1}^{\sigma}/(1+\chi_{m})$
\item $r_{m}= s_{m}-\chi_{m}t_{m}$
\item $r_{m}^{\sigma}=c_{m}^{\sigma}q_{m}^{\sigma}r_{m}$
\item second quasi-minimization and update iterative
\item $\theta_{m}=\|r_{m}^{\sigma}\|/\tilde{\tau}; \ \zeta=1/\sqrt{1+\theta_{m}^{2}}$
\item $\tau=\tilde{\tau}\tilde{\theta_{m}}\zeta; \   \eta_{m}=\zeta^{2}\chi_{m}^{\sigma}$
\item $d_{m}^{\sigma}=s_{m}^{\sigma}+\frac{\tilde{\theta}_{m-1}^{2}\tilde{\eta}_{m-1}}{\chi_{m}^{\sigma}}\tilde{d}_{m-1}$
\item $x_{m}^{\sigma}=\tilde{x}_{m-1}^{\sigma}+\eta_{m}d_{m}$
\item If $x_{m}^{\sigma}$ is accurate enough, then quit
  \end{itemize}
  \item {\bf end for}
\end{enumerate}

\section{ Numerical examples}

In this section, some numerical experiments will be described. The goal of these experiments is to examine the effectiveness of the SQMRCGstab method.

All the numerical experiments were performed in MATLAB 7.1. The machine we have used is a PC-Pentium(R)4, CPU 2.50 GHz, 2.00 GB of RAM. In all of our runs, we used a zero initial guess. All the convergence of the numerical experiments were illustrated in figures. The horizontal axis of figures is the number of matrix-vector multiplies, the vertical axis is relative norm of residual $\parallel r_{m}\parallel/\parallel r_{0}\parallel$.

\textbf{Example 5.1}

Let us consider circuit simulation matrices from Rajat and Raj, which were taken form the University of Florida Sparse Matrix Collection \cite{Davis}. The first matrix is a $1960\times1960$ binary symmetric, and the shift parameter is considered two values, $\sigma=1, 10$. The second matrix is a $1879\times1879$ real unsymmetric, and we consider two values for the shift parameter, $\sigma=0.1, 1$, respectively. Their right-hand side is a unit vector.

\begin{figure}[htbp]
\label{fig1}
\centerline{\includegraphics[width=3.20in,height=2.50in]{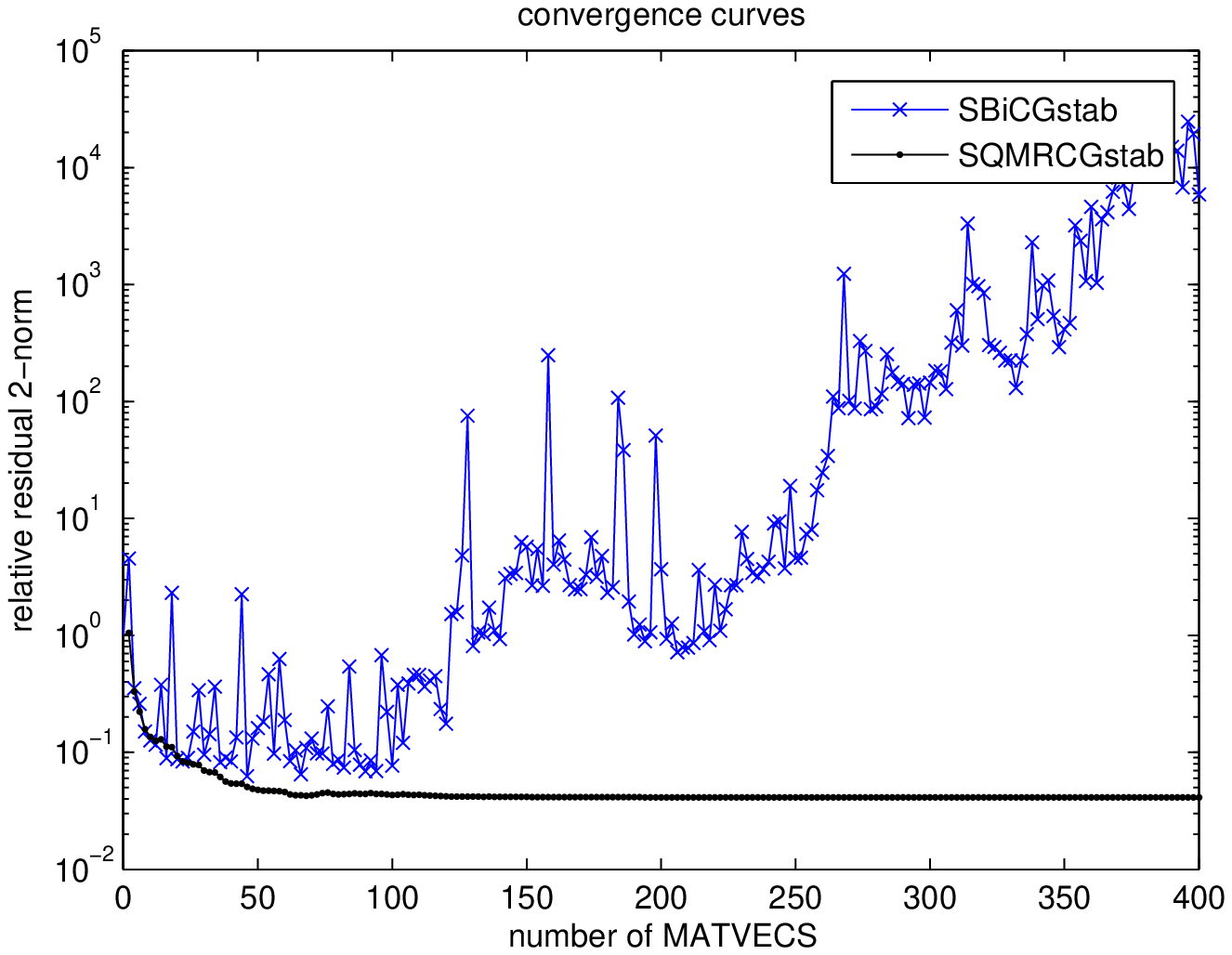}\includegraphics[width=3.20in,height=2.50in]{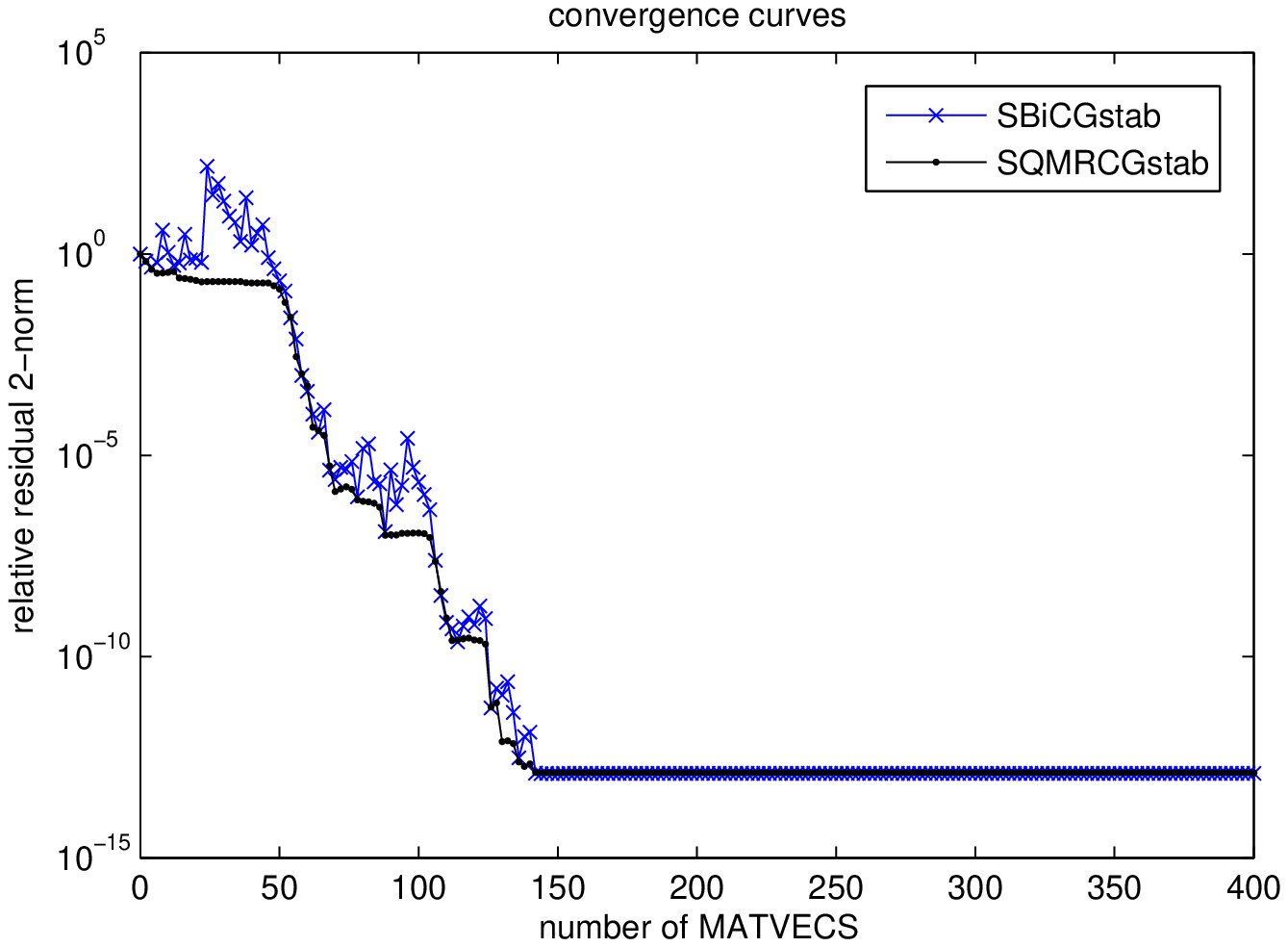}}
\caption{\small Example of binary symmetric matrix. Left: $\sigma=1$, $cond(A+\sigma I)=3.5860e+017$. Right: $\sigma=10$, $cond(A+\sigma I)=60.6806$.}
\end{figure}

\begin{figure}[htbp]
\label{fig1}
\centerline{\includegraphics[width=3.20in,height=2.50in]{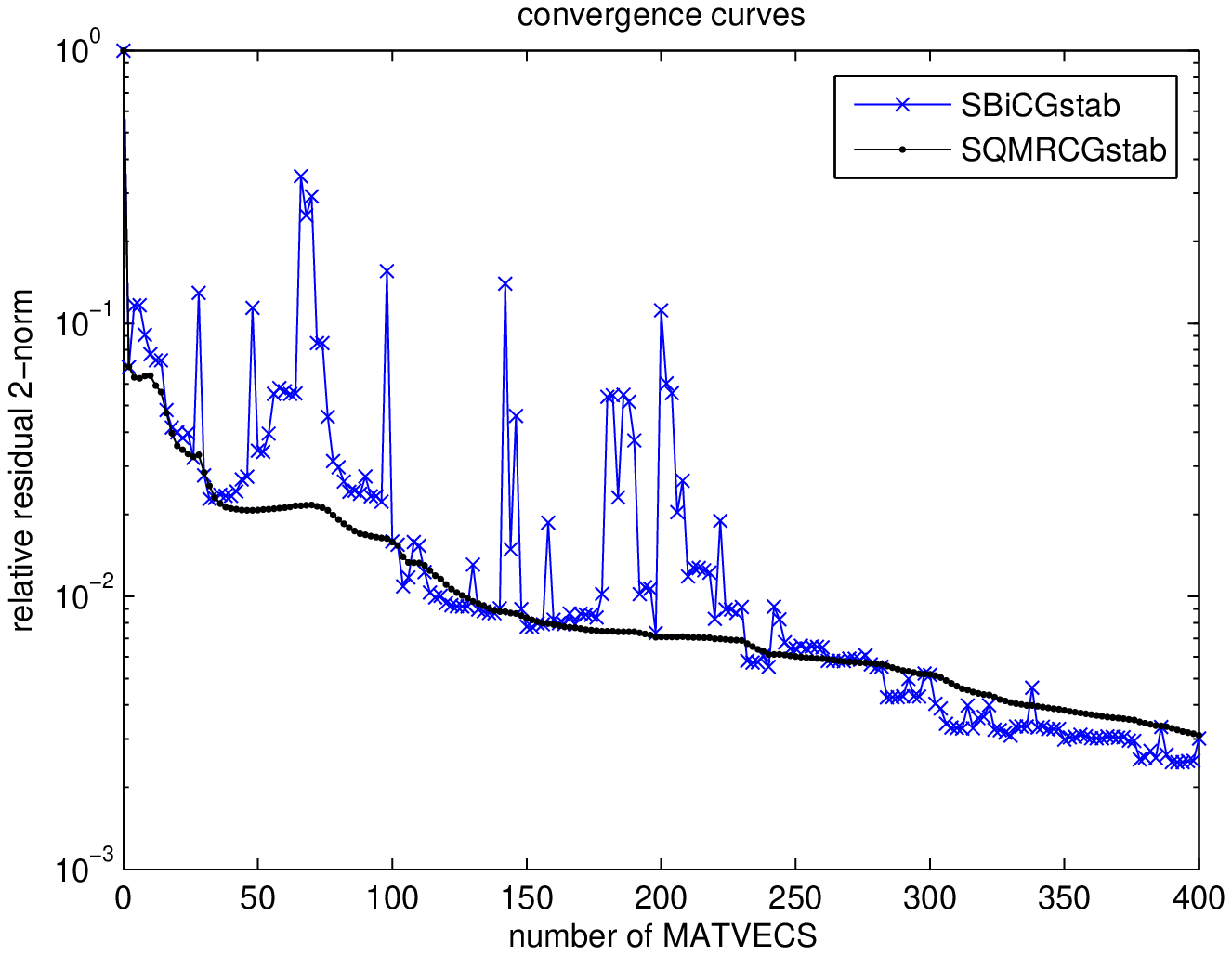}\includegraphics[width=3.20in,height=2.50in]{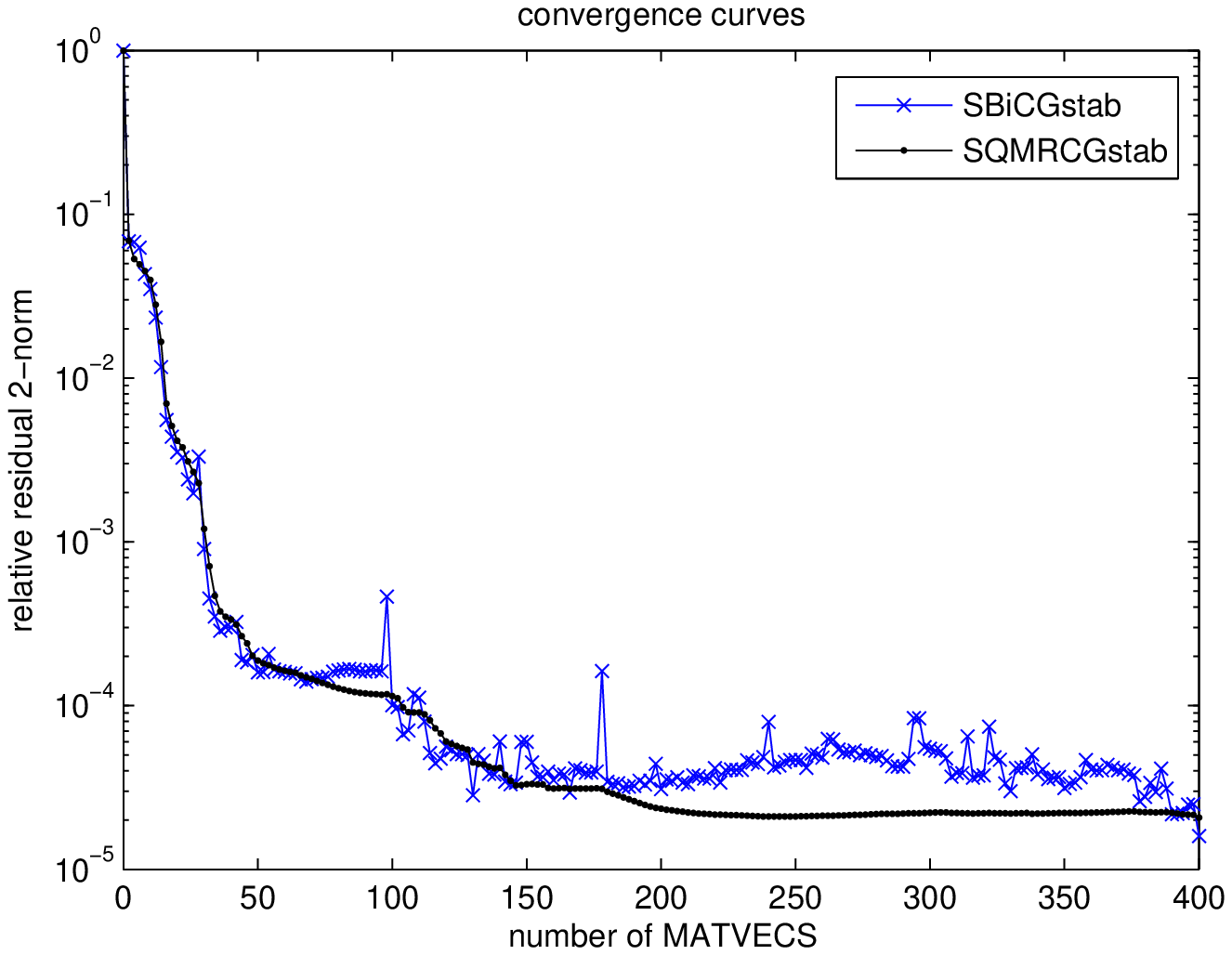}}
\caption{\small Example of real unsymmetric matrix. Left: $\sigma=1$, $cond(A+\sigma I)=3.9112e+004$. Right: $\sigma=10$, $cond(A+\sigma I)=2.1810e+006$.}
\end{figure}

In Figure 1,2, it is observed that the SQMRCGstab method works better than the shifted BiCGstab method, and the convergence plot for the SQMRCGstab method appears well smoothed.  The SQMRCGstab method converges faster even though the condition of the matrix $A+\sigma I$ is more larger (see the left column figures).

\textbf{Example 5.2}

The second numerical experiments stem form a QCD problem. In this part, we compare the SQMRCGstab with the shifted BiCGstab and the Multi-shift QMRIDR(s)\cite{MartinB} methods.

Quark propagators are obtained by solving the inhomogeneous lattice Dirac equation $Ax = b$, where $A = I - kD$ with $0\leq k < k_{c}$ is a large but sparse complex non-Hermitian matrix representing a periodic nearest-neighbour coupling on a four-dimensional Euclidean space-time lattice. The right-hand side vector $b$ is taken as a unite vector.

From the physical theory it is clear that the matrix $A$ should be positive real (all eigenvalues lie in the right half plane) for $0\leq k <k_{c}$. Here, $k_{c}$ represents a critical parameter which depends on the given matrix $D$. It is noted that the matrix $A$ may be a negative real when the parameter $k>k_{c}$.

The matrix and the corresponding critical parameter $k_{c}$ are available from the set QCD of the web repository Matrix Market \cite{Matrix Market}. We take matrices conf5.4-00l4x4-1400.mtx and conf5.4-00l4x4-1800.mtx, which are  $3072\times3072$ complex and non-Hermitian. The corresponding critical values are 0.20328, 0.20265, respectively.  The structures of matrices conf5.4-00l4x4-1400.mtx and conf5.4-00l4x4-1800.mtx are plotted in Figure
3, respectively.

\begin{figure}[htbp]
\label{fig1}
\centerline{\includegraphics[width=3.20in,height=2.50in]{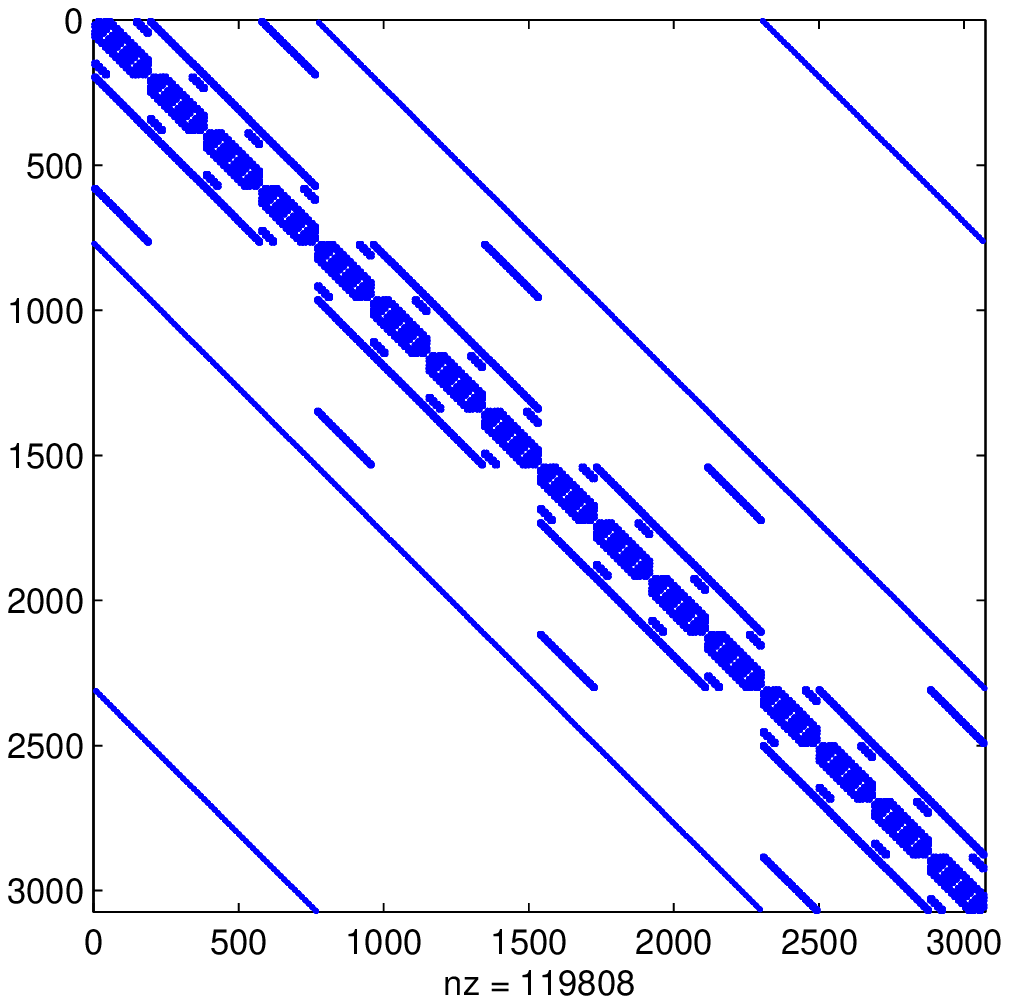}\includegraphics[width=3.20in,height=2.50in]{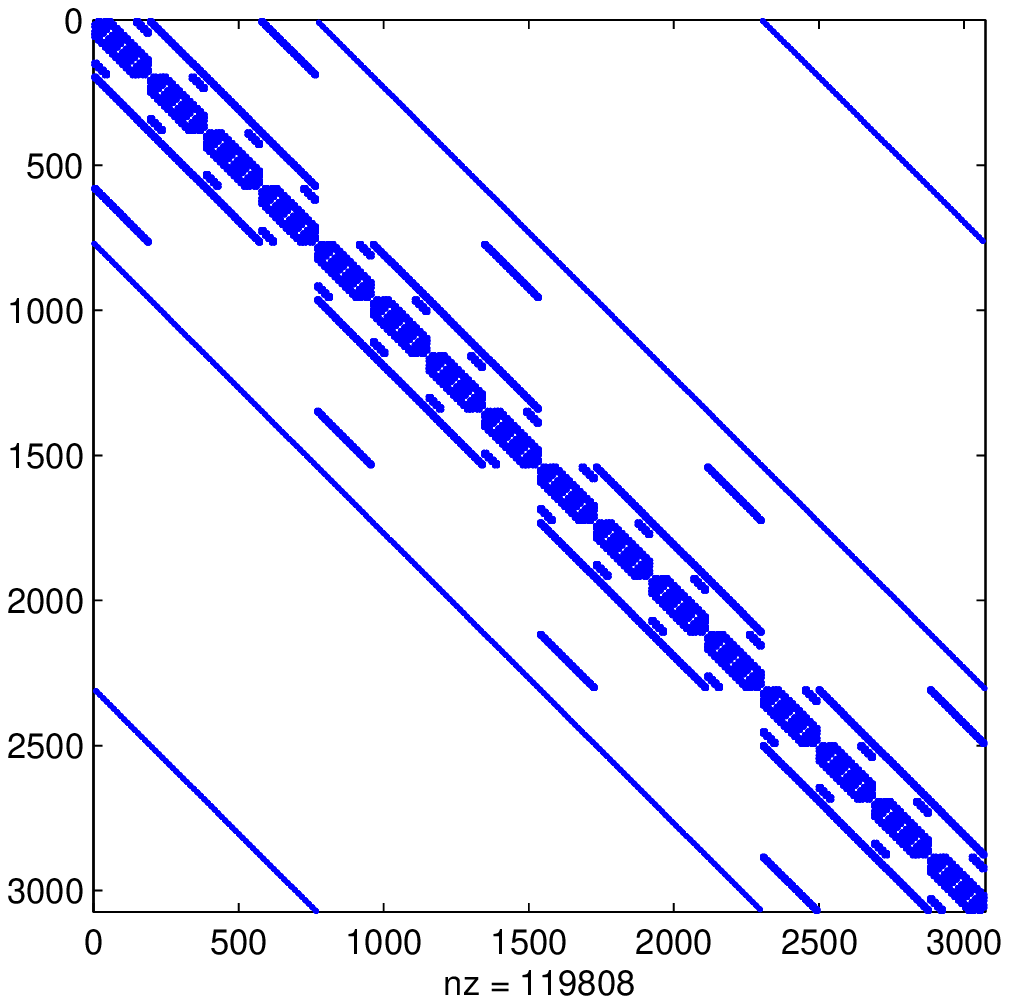}}
\caption{\small Patterns of the matrices in Example 4.2 conf5.4-00l4x4-1400.mtx(left) and
conf5.4-00l4x4-1800.mtx(right).}
\end{figure}

Two different sets of parameters are chosen for numerical experiments. The first set is $k_{1}=0.2$, $k_{2}=0.196$, the second one is $k_{1}=0.2$, $k_{2}=0.176$. In both cases, the seed system is taken to be the system with value $k_{1}$.

As A. Frommer did in \cite{Frommer}, we take a preconditioning process which transforms from the original system to an odd-even-reduced system, which is also a shifted structure. That preconditioning process, which is considered to be the only successful preconditioning in the QCD community so far, usually reduces the number of iterations in a Krylov subspace method. The procedure is described as follows.

First, the grid points are ordered by a red-black (or odd-even) manner, the matrix $D$ becomes
\[
D = \left(
  \begin{array}{cc}   
    0& D_{eo} \\  
   D_{oe} & 0\\  
  \end{array}
\right),                 
\]
Correspondingly,
 \[
x = \left(
  \begin{array}{c}   
    x_{e} \\  
   x_{o}\\  
  \end{array}
\right),  \   \
b = \left(
  \begin{array}{cc}   
    b_{e} \\  
   b_{o}\\  
  \end{array}
\right).                
\]
Second, substituting the above form into this formula $(I - kD)x=b$, the following equation is got
\[
x_{e}-kD_{eo}x_{o}=b_{e},\ \ \
-kD_{oe}x_{e}+x_{o}=b_{o}.
\]
Finally, the odd-even-reduced system is obtained
\[
(I-k^{2}D_{oe}D_{eo})x_{o}= b_{o}+kD_{eo}b_{e}.
\]
Since we worked with the odd-even-reduced system, it means that the matrix is transformed into a $1536\times1536$ complex and non-Hermitian with the corresponding value $k^{2}$, for more details, see \cite{Datta,Frommer1}.

If we write the odd-even-reduced system as in (2), then the shifted parameter $\sigma$ can be taken as $\sigma=k_{2}^{-2}-k_{1}^{-2}>0$.
\begin{figure}[htbp]
\label{fig1}
\centerline{\includegraphics[width=3.20in,height=2.50in]{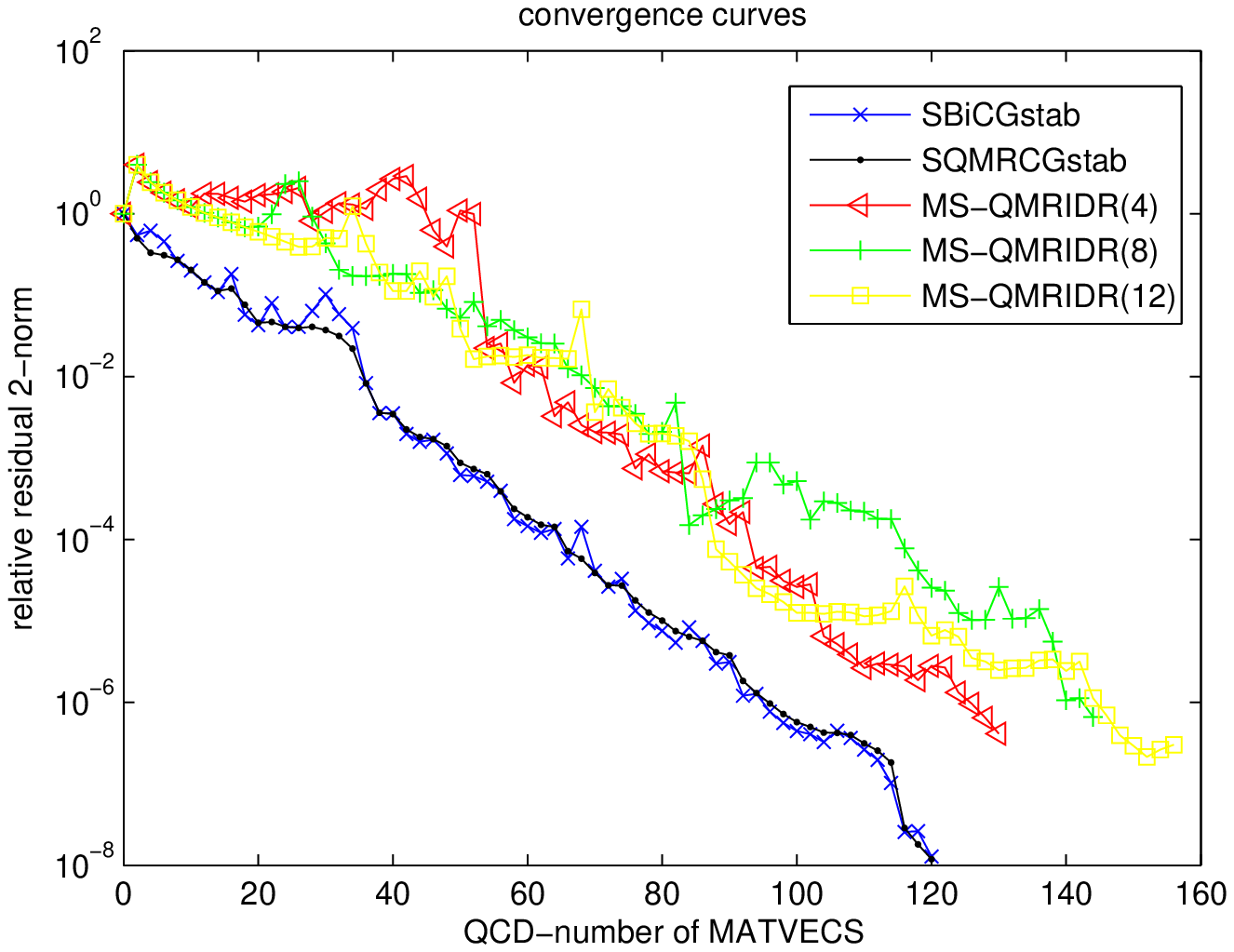}\includegraphics[width=3.20in,height=2.50in]{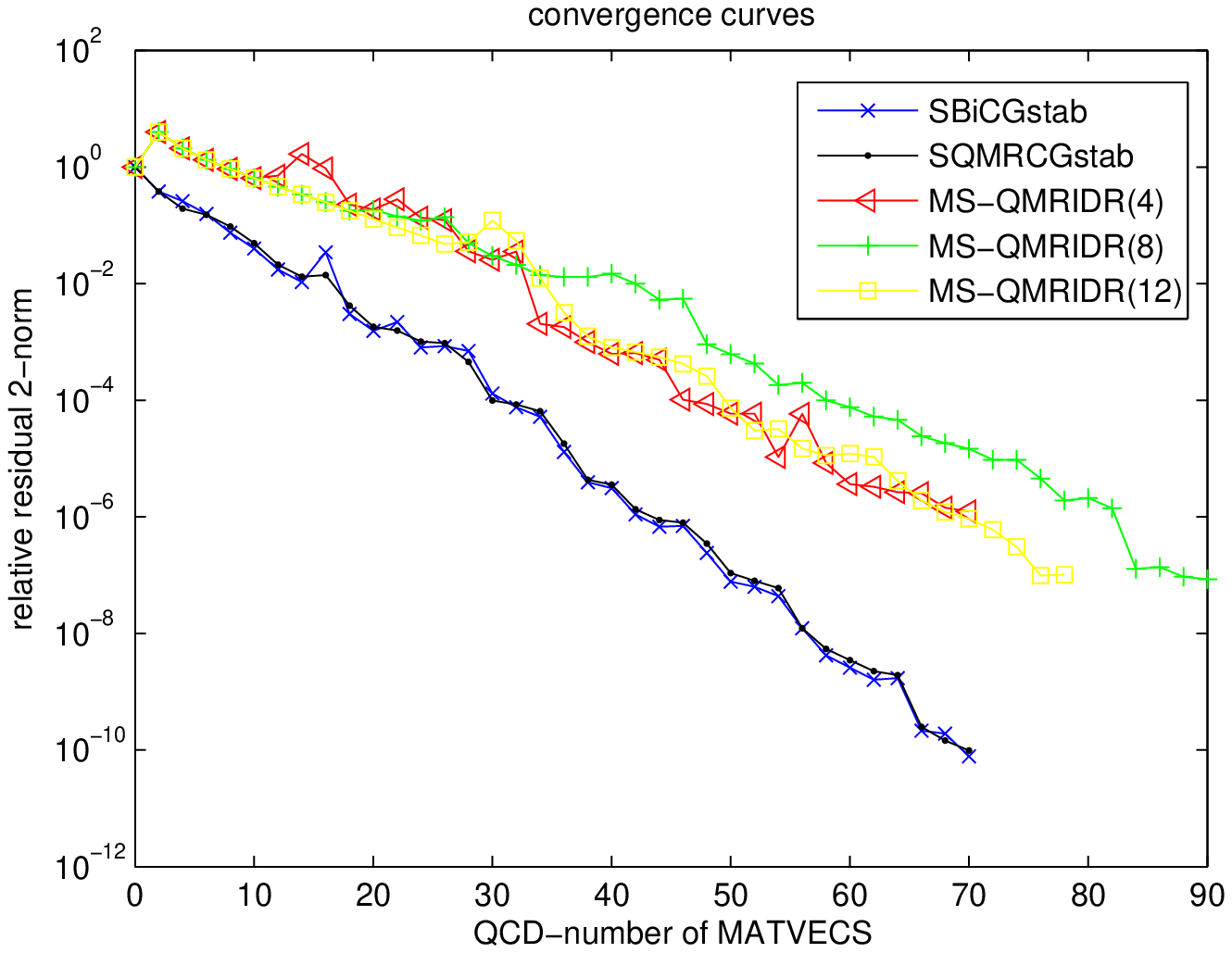}}
\caption{\small conf5.4-00l4x4-1400. Left: $\sigma=0.196^{-2}-0.2^{-2}$, Right: $\sigma=0.176^{-2}-0.2^{-2}$.}
\end{figure}

\begin{figure}[htbp]
\label{fig1}
\centerline{\includegraphics[width=3.20in,height=2.50in]{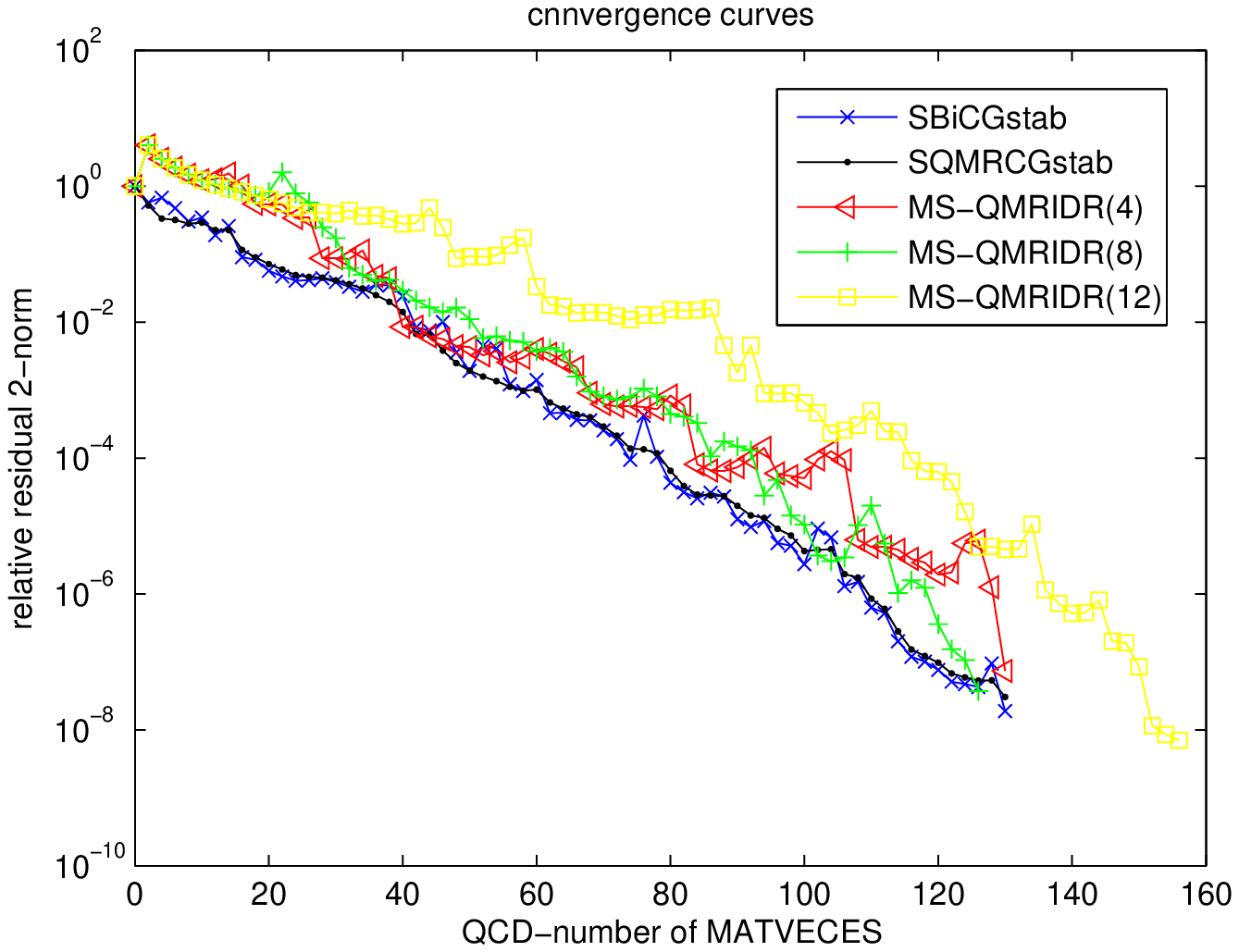}\includegraphics[width=3.20in,height=2.50in]{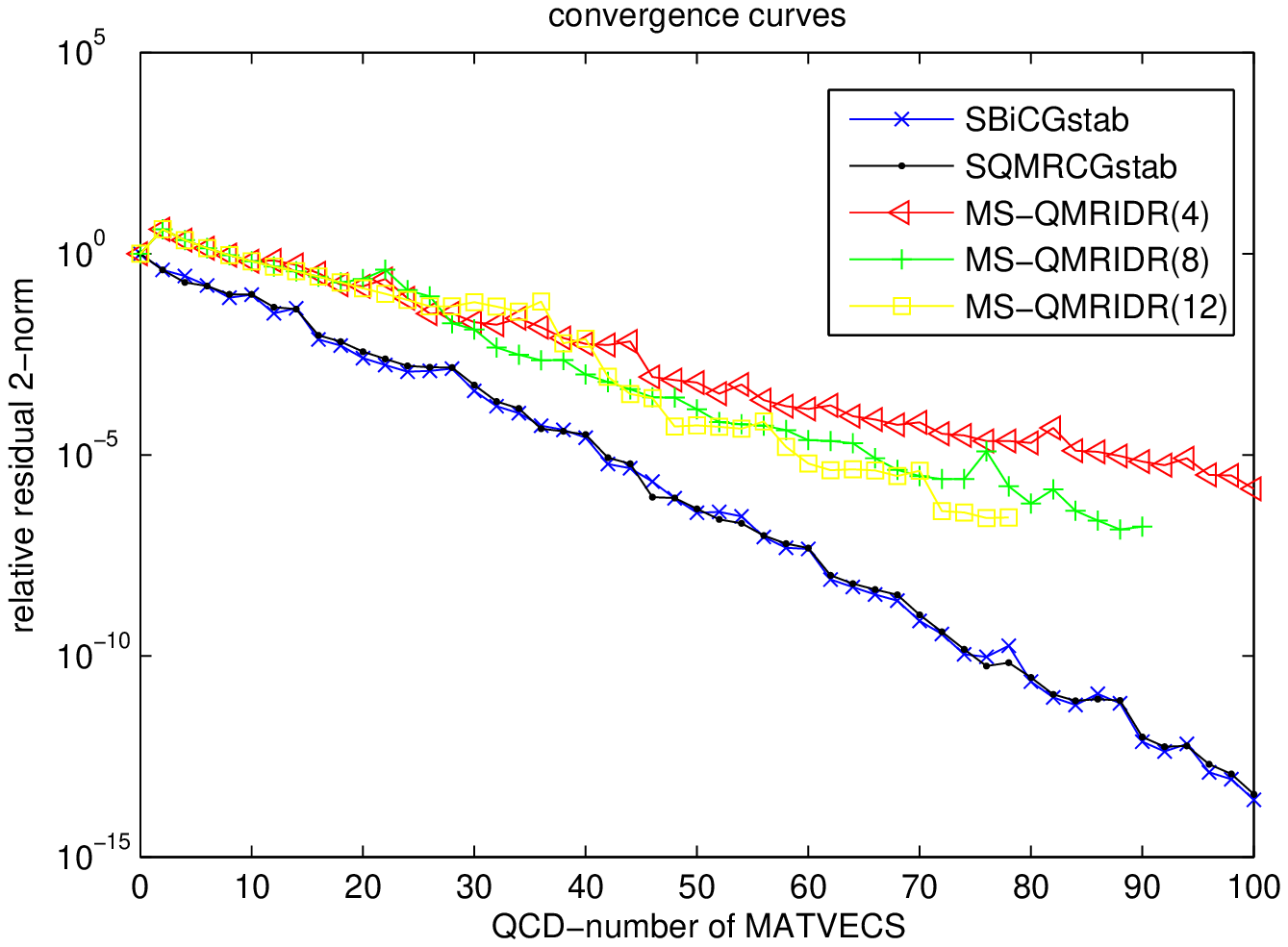}}
\caption{\small conf5.4-00l4x4-1800. Left: $\sigma=0.196^{-2}-0.2^{-2}$, Right: $\sigma=0.176^{-2}-0.2^{-2}$.}
\end{figure}
A. Frommer \cite{Frommer} identified that for a positive real matrix $A$ and a positive shifted parameter $\sigma$,
a damping of the half cycle residual $s^{\sigma}_{m}$ for the shifted system which, albeit not optimal, is larger than the damping obtained on the seed system. Convergence curves of the SQMRCGstab, the shifted BiCGstab and the MS-QMRIDR(s) are displayed in Figure 4, 5. Apparently, the SQMRCGstab and the BiCGstab are more competitive than the MS-QMRIDR(s) on the QCD problem. Moreover, the SQMRCGstab achieves a smoothing of the residual compared to the shifted BiCGstab method.

\textbf{Example 5.3}

The third example stems from the structural dynamics problem \cite{Simoncini}. We consider two cases. The first case is a $100\times100$ upper bidiagonal matrix A with the diagonal the vector $d=[0.001,0.002,0.003,0.004,10+5,11+5,\ldots,105]$ and the super-diagonal the vector of all ones.  The second one is  $1000\times1000$ upper bidiagonal matrix A with diagonal the vector $d=[0.0001,0.0002,0.0003,0.0004,10+5,\ldots,1005]$ and the super-diagonal the vector of all ones. Two values were considered for the shifted parameter, $\sigma=1,-1$, and their right-hand side is the vector of all ones, normalized to have unit norm.

As it can be seen from Figures 6 and 7. By using the SQMRCGstab method, we obtain the more smoother and faster convergent plots.
\begin{figure}[htbp]
\label{fig1}
\centerline{\includegraphics[width=3.20in,height=2.50in]{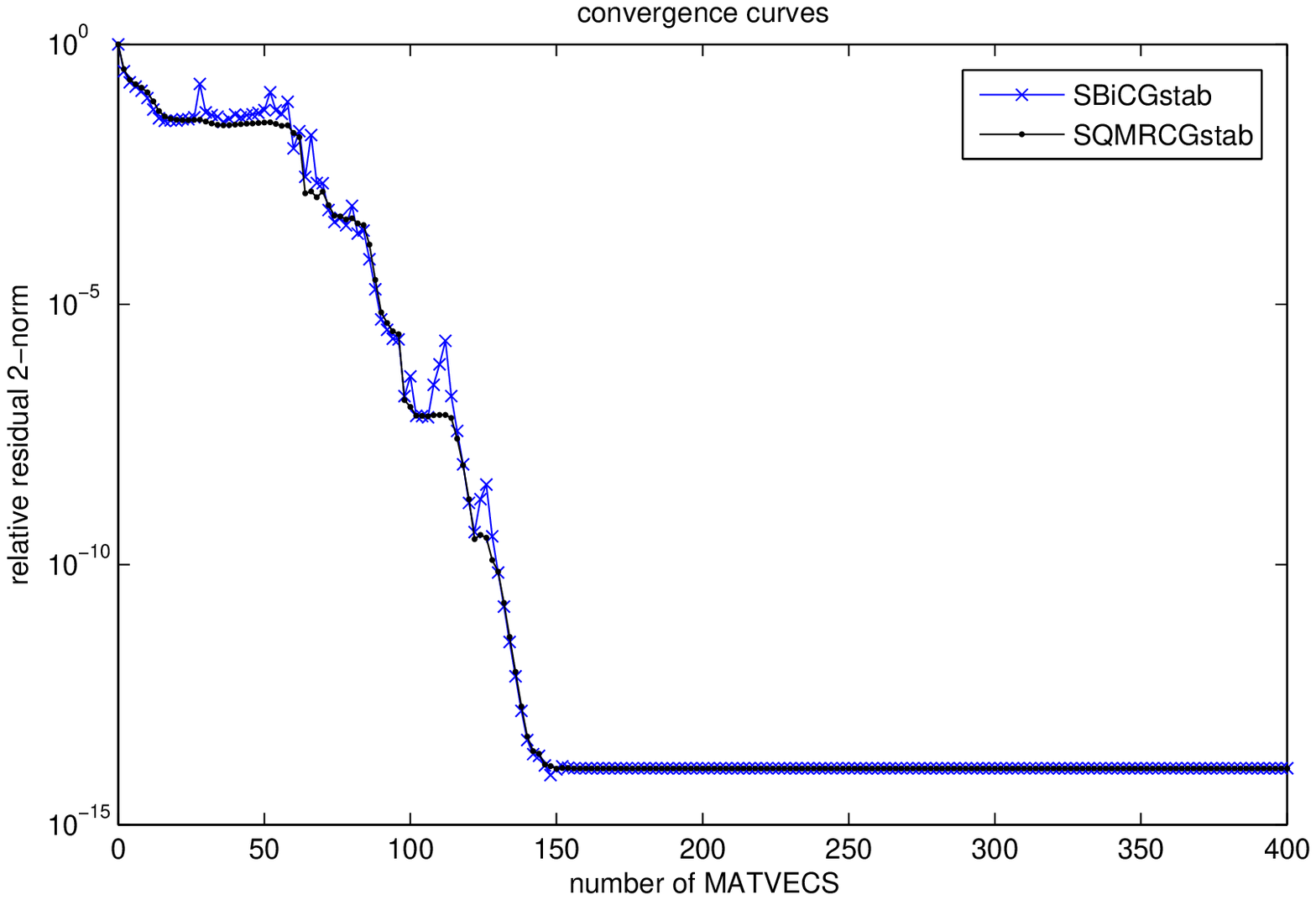}\includegraphics[width=3.20in,height=2.50in]{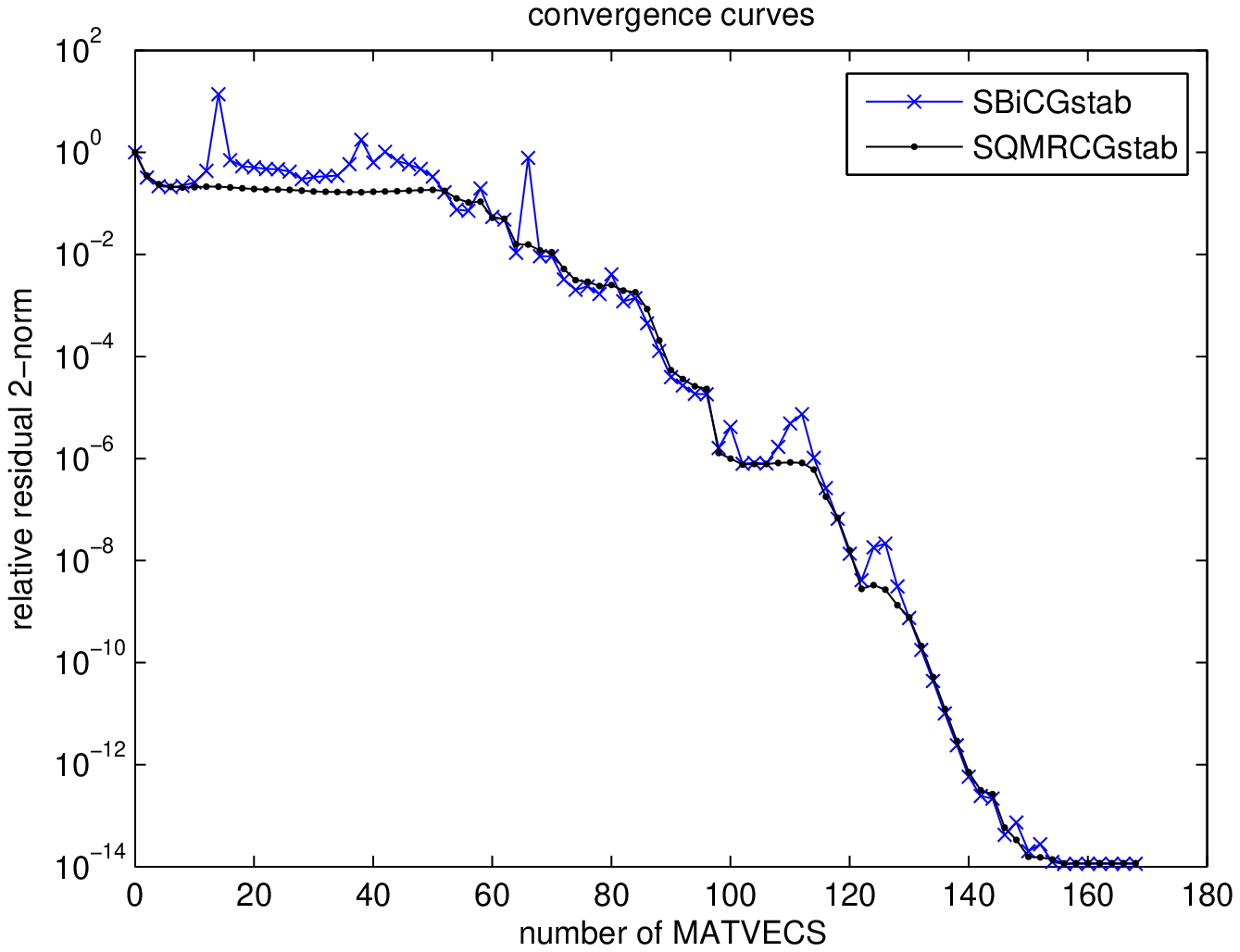}}
\caption{\small Case 1. Left: $\sigma=1$, Right: $\sigma=-1$.}
\end{figure}

\begin{figure}[htbp]
\label{fig1}
\centerline{\includegraphics[width=3.20in,height=2.50in]{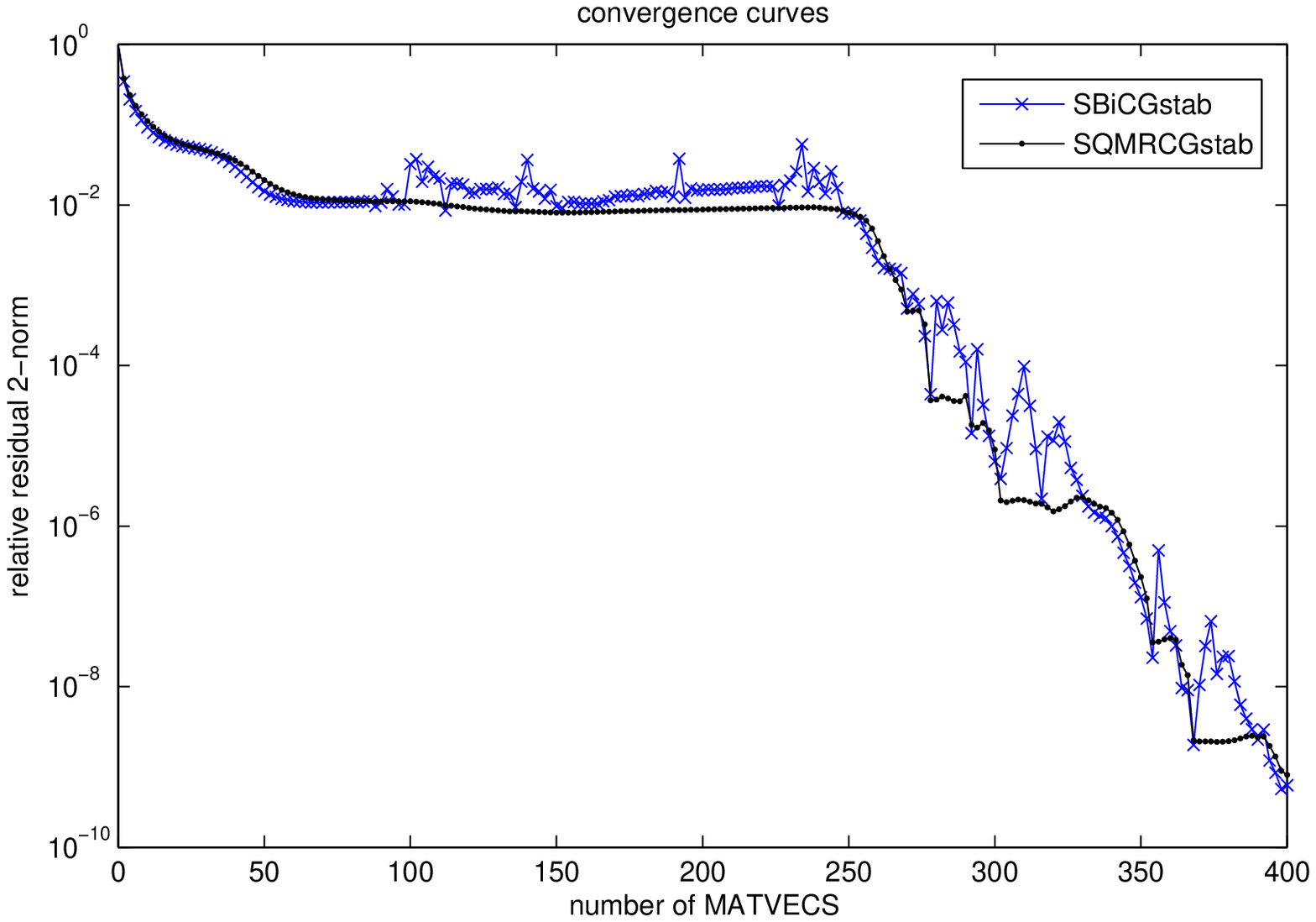}\includegraphics[width=3.20in,height=2.50in]{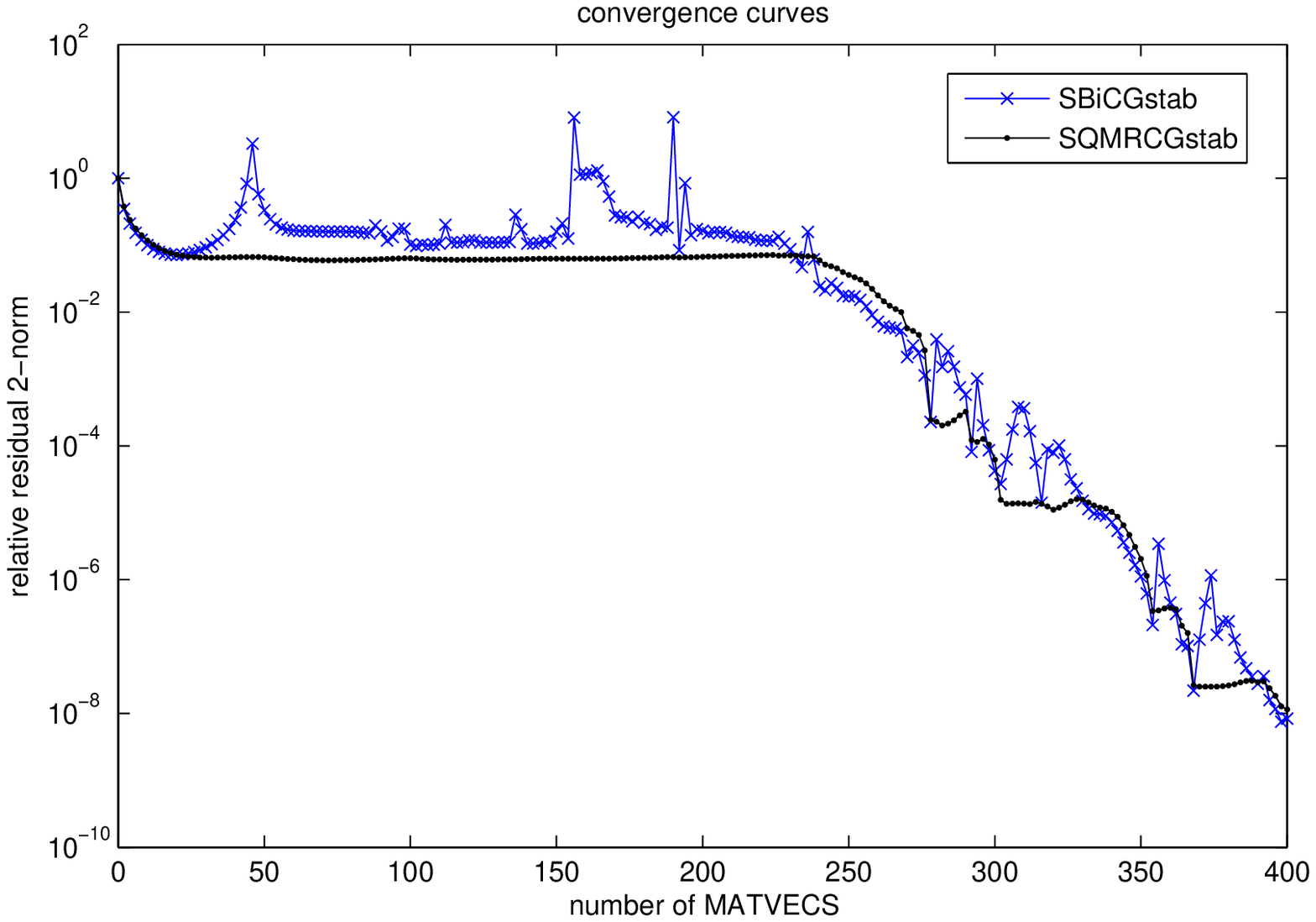}}
\caption{\small Case 2. Left: $\sigma=1$, Right: $\sigma=-1$.}
\end{figure}

\section{Conclusions and future work}

In this paper, we derived a SQMRCGstab method. Our motivation for this method is to inherit any potential improvements on performance BiCGstab,
while at the same time providing a smoother convergence behavior. The SQMRCGstab method has many desirable properties, such as transpose-free, short recurrences. Most important of all, it can make the number of matrix-vector products and the number of inner products be the same as those for a single linear system. Numerical experiments on many real problems confirm the theoretical results and show that our approach is more efficient than the shift BiCGstab method and the MS-QMRIDR(s) method.
The SQMRCGstab method is more competitive than shifted BiCGstab(2) illustrated by some experiments which are not shown up, but inferior to shifted BiCGstab(4) \cite{Frommer}. Therefore, in older to promote the competitiveness of the SQMRCGstab, it lead us to explore several variants of the SQMRCGstab in the future work, just like the shifted BiCGstab($\ell$), which is the shifted BiCGstab's generalization. These problems are important and interest, which will be further investigated and solved in later work.


{\small \begin{thebibliography}{99}

\bibitem{Bloch1} J.C.R. Bloch, T. Wettig, Domain-well and overlap fermions at nonzero quark chemical potential, Phys. Rev. D 76: 114511(2007).

\bibitem{Bloch2} J.C.R. Bloch, T. Breu, A. Frommer, S. Heybrock, etc., Short-recurrence Krylov subspace methods for the overlap Dirac operator at nonzero chemical potential, Comput. Phys. Commun., 181: 1378-1387 (2010).

\bibitem{Bloch3} J.C.R. Bloch, A. Frommer, B. Lang, and T. Wettig, An iterative method to compute the sign function of a non-Hermitian matrix and its application to the overlap Dirac operator at nonzero chemical potential, Comput. Phys. Commun., 177: 933-943 (2007).

\bibitem{Bloch4} J.C.R. Bloch and S. Heybrock, A nested Krylov subspace method to compute the sign function of large complex matrices, Comput. Phys. Commun., 182: 878-889 (2011).

\bibitem{Bloch5} J.C.R. Bloch, S. Heybrock, A nested Krylov subspace method for the overlap operator, proceedings of the XXVII International Symposium on Lattice Field Theory, Bejing, China, PoS(LAT2009) 025 (2009).

\bibitem{Sakurai} T. Sakurai,H. Tadano,Y. Kuramashi, Application of block Krylov subspace algorithms to the Wilson-Dirac equation with multiple right-hand sides in lattice QCD, Comput. Phys. Commun., 181: 113-117 (2010).

\bibitem{Datta} B. Datta, Y. Saad. Arnoldi methods for large Sylvester-like observer matrix equations and an associated algorithm for partial spectrum assignment, Linear Algebra Appl., 154-156:225-244 (1991).

\bibitem{Gallopoulos} E. Gallopoulos, Y. Saad, Efficient parallel solution of parabolic equations, SIAM, Philadelphia,PA,, 251-256 (1990).

\bibitem{Van} J. Van den Eshof, G.L.G. Sleijpen, Accurate conjugate gradient methods for families of shifted systems, Appl. Numer.Math., 49:17-37 (2004).

\bibitem{Takayama} R. Takayama, T. Hoshi, T. Sogabe, Linear algebraic calculation of Green's function for large-scale electronic structure theory, Phys.Rev.B, 73:165108 1-9 (2006).

\bibitem{Sogabe} T. Sogabe, Shao-Liang Zhang, An Extension of the COCR Method to Solving Shifted Linear Systems with Complex Symmetric Matrix, East Asian Journal on Applied Mathematics, 2:97-107 (2011).

\bibitem{Frommer1} A. Frommer, Uwe Gl$\ddot{a}$ssner, Rstarted GMRES For Shifted Linear Systems, SIAM J. Sci. COMPUT., 19:15-26(1998).

\bibitem{Simoncini} V. Simoncini, Restarted Full Orthogonalization Method for Shifted linear systems, BIT Numerical Mathematics, 43:459-466 (2003).

\bibitem{Freund1} R.W Freund, Solution of shifted linear systems by quasi-minimal residual iterations, Numerical Linear Algebra: Proceedings of the Conference in Numerical Linear and Scientific computation Kent (Ohio), New York, (1993).

\bibitem{Jegerlehner} B. Jegerlehner, Krylov space solvers for shifted linear systems, Arxiv preprint hep-lat/9612014, (1996).

\bibitem{Frommer} A. Frommer, BICGStab(l) for Families of Shifted Linear Systems, Computing, 70:87-109 (2003).

\bibitem{Chan} T.F. Chan,E. Gallopoulos, V. Sinoncini, A Quasi-Minimal Residual Variant of the BiCGstab Algorithm for Nonsymmetric Systems, SIAM, 182:81-90 (1994).

\bibitem{Freund} R.T. Freund, A transpose-free quasi-minimal residual algorithm for non-hermitian linear systems, SIAM J. Sc. Stat. Comp., 14:470-482 (1993).

\bibitem{Saad} Y Saad, Iterative methods for sparse linear systems, SIAM,(2003).

\bibitem{Vann} van den Eshof, J., Frommer, A., Lippert, Th., Schilling, K., van der Vorst, H. A.: Numerical methods for the QCD overlap operator I. Sign-function and error bounds, Comp. Physics Comm., 146, 203-224 (2002).

\bibitem{Vorst} van der Vorst, H. A. BI-CGSTAB: A fast and smoothly converging variant of BI-CG for the solution of nonsymmetric linear systems, SIAM J. Sci. Stat. Comput., 13, 631-644 (1992).

\bibitem{Davis} T. A. Davis and Y. F. Hu, The University of Florida Sparse Matrix Collection, ACM Trans. Math. Software, to appear; available online at http://www.cise.ufl.edu/research/sparse/matrices.

\bibitem{Matrix Market} National Institute of Standards and Technology: Matrix Market, http://math.nist.gov/Matrix-Market

\bibitem{Eshof} J. van den Eshof, A. Frommer, T. Lippert, K. Schilling, H.A. van der Vorst , Numerical Methods for the QCD Overlap Operator: I. Sign-Function and Error Bounds, Comp. Physics Comm., 146:203-224 (2002).

\bibitem{Varga} R.S. Varga, Matrix Iterative Analysis. Prentice-Hall, Englewood Cliffs, NJ, 1962.

\bibitem{Arnold} G. Arnold, N. Cundy, J. van den Eshof, A. Frommer, S. Krieg, T. Lippert, K. Sch$\ddot{a}$fer, Numerical Methods for the QCD Overlap Operator: II. Optimal Krylov Subspace Methods, Comput. Appl. Math., 164:587-600 (2004).

\bibitem{Martin} Martin H. Gutknecht, IDR explained, Electron. Trans. Numer. Anal., 36:126-148. (2009/10).

\bibitem{MartinB} Martin B. van Gijzen, Gerard L.G. Sleijpen and Jens-Peter M. Zemke, Flexible and Multi-Shift Induced Dimension Reduction Algorithms for solving Large Sparse Linear Systems. Delft University of Technology, Reports of the Department of Applied Mathematical Analysis, Report 11-06, 2011.

\bibitem{MartinB1} M.B. van Gijzen, P. Sonneveld, An elegant IDR(s) variant that efficiently exploits bi-orthogonality properties, Report 08-21, Department of Applied Mathematical Analysis, Delft University of Technology, 2008.

\bibitem{G.L.G} G.L.G. Sleijpen, M.B. van Gijzen, Exploiting BiCGstab($\ell$) strategies to induce dimension reduction, SIAM J. Sci. Comput., 32 (2010) 2687-2709.
\end{thebibliography}}
\end{document}